\documentclass[11pt]{amsart}
\usepackage{amsmath}
\usepackage[active]{srcltx}
\usepackage{t1enc}
\usepackage[latin2]{inputenc}
\usepackage{verbatim}
\usepackage{amsmath,amsfonts,amssymb,amsthm}
\usepackage[mathcal]{eucal}
\usepackage{enumerate}
\usepackage[centertags]{amsmath}
\usepackage{graphics}

\newtheorem{theorem}{Theorem}

\newtheorem{remark}{Remark}

\newtheorem*{Te}{Theorem T}

\newtheorem*{Ha}{Theorem H}
\newtheorem*{GO}{Theorem GO}
\newtheorem*{KS}{Theorem KS}

\begin{document}
\author{ Ushangi Goginava}
\address{U. Goginava, Department of Mathematical Sciences \\
	United Arab Emirates University, P.O. Box No. 15551\\
	Al Ain, Abu Dhabi, UAE}
\email{zazagoginava@gmail.com; ugoginava@uaeu.ac.ae}
\author{ Farrukh Mukhamedov}
\address{F. Mukhamedov,
	New Uzbekistan University, 54, Mustaqillik ave. \\
	100007, Tashkent, Uzbekistan;
	Central Asian University, Tashkent 111221, Uzbekistan;
	United Arab Emirates University, P.O. Box No. 15551\\
	Al Ain, Abu Dhabi, UAE}
\email{farrukh.m@uaeu.ac.ae}

\title[Divergence of Fourier series]{On problems of the Divergence of logarithmic means of Fourier series}
\date{}

\begin{abstract}
The main objective of the present paper is to solve the several long standing divergence problems related to the logarithmic means of Fourier series in context of general orthonormal systems.
\end{abstract}

\maketitle

\footnotetext{%
}

\section{Introduction}

Let $\psi :[0,\infty )\rightarrow \lbrack 0,\infty )$ be a non-decreasing
function. The class of all measurable functions $f:\mathbb{I}\rightarrow
\mathbb{R}$, $\mathbb{I}:=[0,1)$ for which $\int_{\mathbb{I}}\psi
(|f|)<\infty $ will be denoted by $\psi \left( L\right) \left( \mathbb{I}%
\right) $. If $\psi \left( u\right) =u$ then $\psi \left( L\right) \left(
\mathbb{I}\right)$ coincides with the standard $L_1$-space, i.e.   $\psi \left( L\right) \left(
\mathbb{I}\right):=L_{1}\left( \mathbb{I}\right) $.

In the sequel, by $\Psi $ we denote the set of all non-decreasing functions $\psi
:[0,\infty )\rightarrow \lbrack 0,\infty )$ such that $\liminf\limits_{u%
\rightarrow \infty }\psi (u)/u>0$. Note that for every set $E\subset \mathbb{%
R}$ with finite measure and every function $\psi \in \Psi $, the inclusion $%
\psi (L)(E)\subset L_{1}(E)$ is valid.

It is well-known that Kolmogoroff in \cite{kolmogoroff1923serie} solved the Lusin's problem
by constructing an integrable function $f\in L_{1}(\mathbb{I})$ whose
Fourier trigonometric series diverges almost everywhere. Moreover, in his work \cite{kolmogororf1926s}, Kolmogorov provided an example of a function $f\in L_{1}(\mathbb{I})$ where the divergence occurs everywhere. Similar outcomes were 
demonstrated by Stein for the Walsh-Paley system, as detailed in \cite{stein1961limits}. These results can be summarized as follows.

\begin{KS}
Assume that $S_{n}(f,x)$ is the partial sum of Fourier series of a function $%
f$ with respect to either the trigonometric or the Walsh-Paley system. Then
there is $f\in L_{1}\left( \mathbb{I}\right) $ such that
\begin{equation}
\sup_{n}|S_{n}(f,x)|=\infty ,\ \ \ \text{for all}\ x\in \mathbb{I}.  \label{SS}
\end{equation}
\end{KS}

Let $\Phi={\varphi_n(x)}$ denote an orthonormal system, satisfying $|\varphi_n(x)|\leq M$ for all $x\in \mathbb{I}$ and $n\in \mathbb{N}$. Bochkarev \cite{bochkarev1975fourier} has established the existence of $f\in L_{1}\left( \mathbb{I}\right)$ such that the inequality \eqref{SS} holds for the partial sums of the Fourier series of $f$ concerning the system $\Phi$ on a certain set $E\subset \mathbb{I}$ with $|E|>0$. It is important to note that Kazarian \cite{kazaryan1991problem} demonstrated that the result by Bochkarev does not hold for the set $E$ with $|E|=1$.

The exploration into the divergence of Fourier series remains an ongoing pursuit, with the challenge of identifying the precise class of functions exhibiting divergence. In this context, we highlight certain established findings in this direction. 

In \cite{kon} Konyagin established, for the trigonometric system, the existence of a function $$f\in \psi
(L)\left( \mathbb{I}\right) $$ such that $\sup_{k\in \mathbb{N}}\left\vert
S_{k}(f)(x)\right\vert =\infty $ holds for every $x\in \mathbb{I}$, where $\psi
(u)=o(u\sqrt{\log u/\log \log u})$.

If $\psi (u)=o(u\sqrt{\log u})(u\rightarrow \infty )$, then Bochkarev \cite%
{Bochk} investigated the Walsh-Paley system and established the existence of a function $f\in \psi (L)\left( \mathbb{I}%
\right) $ such that $$\sup_{k\in \mathbb{N}}\left\vert S_{k}(f)(x)\right\vert
=\infty $$ holds for every $x\in \mathbb{I}$.

We recall that the logarithmic means of a function is defined by
\begin{equation}
\label{loga}
L_{n}(f,x)=\frac{1}{\ell_{n}}\sum_{j=0}^{n-1}\frac{S_{j}(f,x)}{n-j},
\end{equation}%
where $\ell_{n}:=\sum\limits_{k=1}^{n}\left( 1/k\right) $.
It is well-known that the convergence of the  logarithmic means is weaker than the convergence of the partial sums.
Consequently, such types of means have been applied in the Fourier series. In 1992, the convergence of logarithmic means of Fourier series was posed as a problem in \cite{MorSiddJAT}. While this problem has been successfully addressed for the Walsh system \cite{GG2006-log,GatGogiAMSDiv2}, the question of almost everywhere convergence/divergence for the trigonometric or general orthonormal system remains open.

The principal objective of this paper is to resolve the aforementioned issues concerning both trigonometric and general orthonormal systems. To achieve this, we have derived a concise and simplified result (refer to Theorem \ref{T4}), enabling us to establish our key findings.
Let us formulate our first result of this paper.

\begin{theorem}
\label{T1} The next statements hold true:

\begin{itemize}
\item[(i)] Let $\psi (u)=o(u\sqrt{\log u/\log \log u})$ then there is $f\in
\psi (L)\left( \mathbb{I}\right) $ such that
\begin{equation}
\sup_{n}|L_{n}(f,x)|=\infty ,\ \ \ \text{for all}\ x\in \mathbb{I},
\label{divv}
\end{equation}%
for the trigonometric system;

\item[(ii)] Let $\psi (u)=o(u\sqrt{\log u})$ then there is $f\in \psi
(L)\left( \mathbb{I}\right) $ such that (\ref{divv}) holds for the
Walsh-Paley system;

\item[(iii)] Let $\Phi =\{\varphi _{n}(x)\}$ be an orthonormal system such
that $|\varphi _{n}(x)|\leq M$ for all $x\in \mathbb{I}$, $n\in \mathbb{N}$.
Then there is $f\in L_{1}\left( \mathbb{I}\right) $ and a set $E\subset
\mathbb{I}$ with $|E|>0$ such that
\begin{equation*}
\sup_{n}|L_{n}(f,x)|=\infty ,\ \ \ \text{for all}\ x\in E.
\end{equation*}
\end{itemize}
\end{theorem}


It is important to emphasize that part (ii) of Theorem \ref{T1} had been previously demonstrated in \cite{GatGogiAMSDiv2}, where the divergence of logarithmic means was established only on a set of positive measure. However, Theorem \ref{T1} not only enhances the previously proven result in \cite{GatGogiAMSDiv2} but also offers a more concise proof for it.

On a different note, Gosselin \cite{gosselin1958divergence} established that for every increasing sequence of natural numbers $\left( n_{k}\right)$, there exists a function $f\in L_{1}(\mathbb{I})$ such that
\begin{equation}
\sup_{k\in \mathbb{N}}\left\vert S_{n_{k}}(f,x)\right\vert =\infty
\label{SZ}
\end{equation}
for almost every $x\in \mathbb{I}$. Conversely, Totik \cite{totik1982divergence} found a function $f\in L_{1}(\mathbb{I})$ for which \eqref{SZ} holds for every $x\in \mathbb{I}$, considering any increasing sequence of natural numbers $\left( n_{k}\right)$.

It turns out that establishing an analogue of the results by Gosselin and Totik remains an open problem, intricately linked with Konyagin's problem \cite{konyagin2006almost}. While a partial solution has been provided in \cite{goginava2021divergence}. 
Let us recall some necessary notions to formulate it.

The \emph{binary coefficients }of a number $n\in \mathbb{N}$ will be denoted
by $\varepsilon _{j}(n)$ $(j\in \mathbb{N}_0)$, i.e., $\varepsilon _{j}(n)\in
\{0,1\}$ for every $j\in \mathbb{N}_{0}$ and $n=\sum_{j=0}^{\infty
}\varepsilon _{j}(n)2^{j}$. Note that $\varepsilon _{j}(n)=0$ if $j$ is
sufficiently large. For $n\in \mathbb{N}$, let us denote
\begin{equation*}
V(n)=\varepsilon _{0}(n)+\sum_{j=1}^{\infty }|\varepsilon
_{j}(n)-\varepsilon _{j-1}(n)|.
\end{equation*}%
The quantity $V(n)$ is called the \emph{variation} of the number $n$. \

Let us define the \emph{spectrum} of a number $n\in \mathbb{N}$ as follows
\begin{equation*}
\text{Sp}(n)=\{j\in \mathbb{N}:\varepsilon _{j}(n)=1\}.
\end{equation*}%
Let $(n_{k})$ be an increasing sequences of natural numbers. We will call $%
(n_{k})$ a sequence with \emph{nested spectrums} if
\begin{equation*}
\text{Sp}(n_{k+1})\cap \lbrack 0,\max \text{Sp}(n_{k})]=\text{Sp}%
(n_{k})\;\;\;\text{for every}\;\;\;k\in \mathbb{N}.
\end{equation*}

\begin{GO}\cite{goginava2021divergence}
\label{t:1} Let $(n_{k})$ be a sequence of natural numbers with nested
spectrums having unbounded variation. Then there exists a function $f\in
L_{1}\left( \mathbb{I}\right) $ such that $$\sup_{k\in \mathbb{N}%
}|S_{n_{k}}(f,x)|=\infty $$ for every $x\in \mathbb{I}.$
\end{GO}

Next our result improves the above mentioned Theorems for the logarithmic means.

\begin{theorem}
\label{T2} The next statements hold true:

\begin{itemize}
\item[(i)] for every increasing sequence of natural numbers $\left(
n_{k}\right) $ there exists a function $f\in L_{1}(\mathbb{I})$ such that
\begin{equation}
\sup_{N}\bigg|\frac{1}{\ln N}\sum_{k=0}^{N-1}\frac{S_{n_{k}}(f,x)}{N-k}\bigg|%
=\infty ,\ \ \ \text{for all}\ x\in \mathbb{I}.  \label{SZ1}
\end{equation}%
for the trigonometric system;

\item[(ii)] Let $(n_{k})$ be nested sequence which has unbounded variation.
Then there is $f\in L_{1}(\mathbb{I})$ such that \eqref{SZ1} holds for the
Walsh-Paley system.
\end{itemize}
\end{theorem}


\begin{remark}
Recently, G\'at \cite{Gat2023} has provided necessary
and sufficient conditions on subsequence $(n_{k})$, for the trigonometric system, such that
\begin{equation}\label{GG1}
\frac{1}{N}\sum_{k=0}^{N-1}S_{n_{k}}(f,x)\rightarrow f\ \ \text{a.e.}
\end{equation}%
for every $f\in L_{1}(\mathbb{I})$. However, from Theorem \ref{T1}  we infer
that the arithmetic means in \eqref{GG1} cannot be replaced by the logarithmic ones.
\end{remark}

Now, let us consider any lacunary sequence $(n_k)$. From Theorem \ref{T2}
one can ask: find an optimal class of functions for which the sequence
\begin{equation*}
\frac{1}{\ln N}\sum_{k=0}^{N-1}\frac{S_{n_k}(f,x)}{N-k}
\end{equation*}
converges a.e. to $f$. Our result allows us to provide an affirmative answer
to this question. Namely, we can prove the next result.

\begin{theorem}
\label{T3} Let $(n_{k})$ be an arbitrary lacunary sequence. Then the space \\$%
L\log ^{+}\log ^{+}L\log ^{+}\log ^{+}\log ^{+}\log ^{+}L\left( \mathbb{I}%
\right) $ is the optimal Orlicz space for a.e. convergence of
\begin{equation*}
\frac{1}{l_{N}}\sum_{k=0}^{N-1}\frac{S_{n_{k}}(f,x)}{N-k}\rightarrow f\text{
\ as }n\rightarrow \infty .
\end{equation*}
\end{theorem}

The proof follows immediately from Theorem \ref{T4} and the results of Di
Plinio \cite{di2014lacunary} and Lie \cite{lie2012pointwise}.

On the other hand, Buzdalin \cite{Buzdalin} provided an affirmative response to Ul'yanov's question by demonstrating that for any set $\mathrm{E} \subset \mathbb{I}$ with $|\mathrm{E}|=0$, there exists a real continuous function $f \in \mathrm{C}(\mathbb{I})$ such that its Fourier series with respect to the trigonometric system diverges unboundedly on $E$. Now, leveraging Theorem \ref{T4}, we can derive our subsequent result.

 \begin{theorem}
  For any set $\mathrm{E} \subset\mathbb{I}$ with $|\mathrm{E}|=0$ there is a real function $f \in \mathrm{C}(\mathbb{I})$ such that its logarithmic means (\ref{loga}) of Fourier series with respect to trigonometric system  diverges unboundedly on $E$.
  \end{theorem}

  Let us point out that an analogue of Buzdalin's theorem remains an open problem for the Walsh systems. Conversely, for bounded functions, this problem has been resolved by Bughadze \cite{Buga}, and the divergence of arithmetic means on sets of measure zero has been explored in \cite{gogimeasurezero}. Karagulyan \cite{KaragulyanMeasureZero} investigated a more general sequence of operators and studied its divergence for bounded functions on a measure zero set.

It is noteworthy that the almost everywhere convergence and divergence of the partial sums of Fourier series have been investigated by many authors, including Setchkin \cite{stechkin1951convergence}, Katznelson \cite{katznelson1966ensembles}, Kahane and Katznelson \cite{kahankatznelson1966ensembles}, Konyagin \cite{Kon2011,Kon2006,Kon2005}, Do and Lacey \cite{DoLacey2012}, and G\'at \cite{gatJAT2010,gatConstr2019}.

\section{Proof of  Main Result}

Let $\mathbf{S}:=\left\{ s_{n}:n\in \mathbb{N}\right\} $ be a sequence of
real numbers. The logarithmic means of the sequence $\mathbf{S}$ are defined
as follows%
\begin{equation}
L_{n}\left( \mathbf{S}\right) :=\frac{1}{\ell_{n}}\sum\limits_{k=0}^{n-1}\frac{%
s_{k}}{n-k}.  \label{log}
\end{equation}
where, as before, $\ell_{n}:=\sum\limits_{k=1}^{n}\left( 1/k\right) $.

Now, let us consider more general means. Assume that $\mathbf{T}:=\left( t_{k,n}\right) $ is an infinite matrix. Define
the $n$th matrix transform of the sequence $\mathbf{S}$ by
\begin{equation}
T_{n}\left( \mathbf{S}\right) :=\sum_{k=0}^{n}t_{k,n}s_{k}\quad (n\in {%
\mathbb{N}}).  \label{T}
\end{equation}%
We say that $T$ is \textit{regular} if $T_{n}\left( \mathbf{S}\right) \rightarrow
s\left( n\rightarrow \infty \right) $ whenever $s_{n}\rightarrow s\left(
n\rightarrow \infty \right) .$ In order that $T$ be regular it is neccessary
and sufficient that (see \cite[Ch. 3]{zygmund2002trigonometric} or \cite[P.
43-57 ]{hardy2000divergent}):

\begin{itemize}
\item[(a)] $\sum_{k=0}^{n}t_{k,n}=1;$

\item[(b)] $t_{k,n}\rightarrow 0$ for each $k$, when $n\rightarrow \infty ;$

\item[(c)] $\sup\limits_{n}\sum\limits_{k=0}^{n}\left\vert
t_{k,n}\right\vert \leq C<\infty $. 
\end{itemize}

We notice that if $t_{n,k}=q_{n-k}/Q_n$, where $(q_n)$ is some sequence,
and
$Q_{n}:=\sum\limits_{k=0}^{n}q_{k}$, then
$T_{n}\left( \mathbf{S}\right)$ coincidences with the N\"{o}rlund means, i.e.
\begin{equation}
	T_{n}\left( \mathbf{S}\right)=: N_{n}^{\left( q\right) }\left( \mathbf{S}\right)=\frac{1}{Q_{n}}%
	\sum\limits_{k=0}^{n}q_{n-k}s_{k},  \label{Nor}
\end{equation}%

If $q_n>0$ for all $n\geq 0$, then the regularity of the N\"{o}rlund means is equivalent to
\begin{equation}\label{qQ}
\lim_{n\to\infty}\frac{q_{n}}{Q_n}=0.
\end{equation}


The $%
\left( C,\alpha _{n}\right) $ means (see  \cite{Kap} for more details) of the sequence $\mathbf{S}:=\left\{
s_{n}:n\in \mathbb{P}\right\} $ is given by

\begin{equation*}
\sigma _{n}^{\alpha _{n}}(\mathbf{S})=\frac{1}{A_{n}^{\alpha _{n}}}%
\sum\limits_{j=0}^{n}A_{n-j}^{\alpha _{n}-1}s_{j},
\end{equation*}%
where
\begin{equation*}
A_{n}^{\alpha _{n}}:=\frac{(1+\alpha _{n})\dots (n+\alpha _{n})}{n!}.
\end{equation*}%
The connections between method $\left( C,\alpha _{n}\right) $ with varying
parameters and Cesaro's method $\left( C,\alpha \right) $ with constant
parameters were studied by Kaplan \cite{Kap}. Regarding, the connection between two
different methods $\left( C,\alpha _{n}\right) $ and $\left( C,\beta
_{n}\right) $ was established recently by Akhobadze and Zviadadze \cite{AkhZviad2019}.

In the sequel, we are going to employ next well-known result (see \cite[pp. 68-69]{hardy2000divergent}).

\begin{Ha} Let $(q_n)$ be a sequence such that
\begin{equation}
q_{0}=1,q_{n}>0,\frac{q_{n+1}}{q_{n}}\geq \frac{q_{n}}{q_{n-1}}\left(
n>0\right) ,  \label{cond0}
\end{equation}%
Assume that \eqref{qQ} holds. Denote
$$
q(x)=\sum_{n=0}^\infty q_nx^n.
$$
Then
\begin{equation}
\frac{1}{q\left( x\right) }=\sum\limits_{n=0}^{\infty }\gamma _{n}x^{n}
\label{=1}, \ \ \ |x|<1,
\end{equation}%
where $\gamma _{0}=1$, $\gamma _{n}<0\left( n>0\right),$ and
\begin{equation}
\sum\limits_{n=1}^{\infty }\gamma _{n}\geq -1.
\label{<1}
\end{equation}
\end{Ha}

Using the previous result we are going to prove next auxiliary fact.

\begin{theorem}
	\label{2met}
Let $\mathbf{S}:=\left\{ s_{n}:n\in \mathbb{P}\right\} $ be a sequence
of real numbers such that 
\begin{equation*}
L_{n}\left( \mathbf{S}\right) \rightarrow s\text{ \ as \ }n\rightarrow
\infty .
\end{equation*}
Then for $\alpha _{n}\in \left( 0,1\right) $
\begin{equation*}
\sigma _{n}^{\alpha _{n}}(\mathbf{S})\rightarrow s\text{ \ \ as \ }%
n\rightarrow \infty .
\end{equation*}%
\end{theorem}
	.
\begin{proof} Let $m$ be any fixed number. We first notice that the condition \eqref{qQ}  (see \cite{hardy2000divergent}) implies the convergence of
	the series $\sum\limits_{n=0}^{\infty }s_{n}x^{n}$ when $\left\vert
	x\right\vert <1$.  Moreover, from $\lim\limits_{n\rightarrow \infty }\frac{A_{n}^{\alpha _{m}-1}}{%
		A_{n-1}^{\alpha _{m}-1}}=1$ we infer the convergence of the series $%
	\sum\limits_{n=0}^{\infty }A_{n}^{\alpha _{m}-1}x^{n}$ when $\left\vert
	x\right\vert <1$. Hence, by the Mertens Theorem one concludes that the next power series is convergent (when  $|x|<1$)
	\begin{eqnarray}\label{rep}
\left( \sum\limits_{n=0}^{\infty }A_{n}^{\alpha _{m}-1}x^{n}\right)
\left( \sum\limits_{n=0}^{\infty }s_{n}x^{n}\right)&=&
\sum\limits_{n=0}^{\infty }\left( \sum\limits_{k=0}^{n}A_{n-k}^{\alpha
	_{m}-1}s_{k}\right) x^{n}\nonumber\\
&=&\sum\limits_{n=0}^{\infty }A_{n}^{\alpha _{m}}\sigma _{n}^{\alpha _{m}}(%
\mathbf{S})x^{n} .
\end{eqnarray}%
Let the sequence $\left( q_{n}\right)
_{n=1}^{\infty }$ satisfy the conditions of Theorem H. Then by Theorem H together with (\ref{rep}) we obtain
\begin{eqnarray*}
\sum\limits_{n=0}^{\infty }A_{n}^{\alpha _{m}}\sigma _{n}^{\alpha _{m}}(%
\mathbf{S})x^{n}
&=&\frac{\left( \sum\limits_{n=0}^{\infty }A_{n}^{\alpha _{m}-1}x^{n}\right)
}{\left( \sum\limits_{n=0}^{\infty }q_{n}x^{n}\right) }\left(
\sum\limits_{n=0}^{\infty }q_{n}x^{n}\right) \left(
\sum\limits_{n=0}^{\infty }s_{n}x^{n}\right) \\
&=&\left( \sum\limits_{n=0}^{\infty }A_{n}^{\alpha _{m}-1}x^{n}\right)
\left( \sum\limits_{n=0}^{\infty }\gamma _{n}x^{n}\right)
\sum\limits_{n=0}^{\infty }\left( \sum\limits_{k=0}^{n}q_{n-k}s_{k}\right)
x^{n} \\
&=&\sum\limits_{n=0}^{\infty }\left( \sum\limits_{k=0}^{n}A_{n-k}^{\alpha
_{m}-1}\gamma _{k}\right) x^{n}\sum\limits_{n=0}^{\infty }Q_{n}N_{n}^{\left(
q\right) }\left( \mathbf{S}\right) x^{n}.
\end{eqnarray*}%
Denote
\begin{equation*}
b_{n,m}:=\sum\limits_{k=0}^{n}A_{n-k}^{\alpha _{m}-1}\gamma _{k}.
\end{equation*}%
From (\ref{rep}) it follows that 
\begin{eqnarray}
\sum\limits_{n=0}^{\infty }A_{n}^{\alpha _{m}}\sigma _{n}^{\alpha _{m}}(%
\mathbf{S})x^{n} &=&\sum\limits_{n=0}^{\infty
}b_{n,m}x^{n}\sum\limits_{n=0}^{\infty }Q_{n}N_{n}^{\left( q\right) }\left(
\mathbf{S}\right) x^{n}  \label{rep2} \\
&=&\sum\limits_{n=0}^{\infty }\left(
\sum\limits_{k=0}^{n}b_{n-k,m}Q_{k}N_{k}^{\left( q\right) }\left( \mathbf{S}%
\right) \right) x^{n}.  \notag
\end{eqnarray}%
Now, we may assume that $n=m$. Then by (\ref{rep2}) one finds
\begin{equation*}
\sigma _{n}^{\alpha _{n}}(\mathbf{S})=\sum\limits_{k=0}^{n}\frac{%
b_{n-k,n}Q_{k}}{A_{n}^{\alpha _{n}}}N_{k}^{\left( q\right) }\left( \mathbf{S}%
\right) .
\end{equation*}

Therefore, the sequence $\left( \sigma _{n}^{\alpha _{n}}(\mathbf{S})\right)
_{n=1}^{\infty }$ can be represented as a matrix transformation $\left(
t_{k,n}:=\frac{b_{n-k,n}Q_{k}}{A_{n}^{\alpha _{n}}}\right) $ of the sequence
$\left( N_{k}^{\left( q\right) }\left( \mathbf{S}\right) \right)
_{n=1}^{\infty }$.

Let us establish the regularity of this transformation. It is enough to prove that conditions $\left( a\right) -\left(
c\right) $ hold for the sequence $t_{k,n}$. Indeed, assume that $%
\mathbf{S}$ is an identical sequence $s_{n}=1,n\in \mathbb{N}$. $\ $Since $%
\sigma _{n}^{\alpha _{n}}(\mathbf{S})=N_{k}^{\left( q\right) }\left( \mathbf{%
S}\right) =1$ we get $\sum\limits_{k=0}^{n}t_{k,n}=1.$ Hence, condition $(a)$
is satisfied. Now, we prove the condition $\left( b\right) $. For fixed $k$ we can
write%
\begin{equation*}
t_{k,n}=\frac{A_{n-k}^{\alpha _{n}}Q_{k}}{A_{n}^{\alpha _{n}}}\frac{1}{%
A_{n-k}^{\alpha _{n}}}\sum\limits_{l=0}^{n-k}A_{n-k-l}^{\alpha _{n}-1}\gamma
_{l}.
\end{equation*}%
Due to  $\frac{A_{n-k}^{\alpha _{n}}Q_{k}}{A_{n}^{\alpha _{n}}}\leq C_{k}$,$\ $%
and the regularity of the method $\left( C,\alpha _{n}\right) $, for  $\gamma
_{k}\rightarrow 0\left( k\rightarrow \infty \right) $ (see (\ref{<1})) we infer that
\begin{equation*}
\frac{1}{A_{n-k}^{\alpha _{n}}}\sum\limits_{l=0}^{n-k}A_{n-k-l}^{\alpha
_{n}-1}\gamma _{l}\rightarrow 0
\end{equation*}%
as $n\rightarrow \infty $. Hence, the condition $\left( b\right) $ holds. It
remains to prove the condition $(c)$.

Assume that the sequence
$\left( q_{n}\right) _{n=1}^{\infty }$ satisfies the following additional
condition%
\begin{equation}
\frac{q_{n}}{q_{n-1}}\leq \frac{A_{n}^{\alpha _{n}-1}}{A_{n-1}^{\alpha
_{n}-1}}=1-\frac{1-\alpha _{n}}{n}\left( n\in \mathbb{N}\right) .
\label{cond}
\end{equation}%
From (\ref{=1}) it follows that
\begin{eqnarray*}
1 &=&\left( \sum\limits_{n=0}^{\infty }\gamma _{n}x^{n}\right) \left(
\sum\limits_{n=0}^{\infty }q_{n}x^{n}\right)  \\
&=&\sum\limits_{n=0}^{\infty }\left( \sum\limits_{k=0}^{n}q_{n-k}\gamma
_{k}\right) x^{n}.
\end{eqnarray*}%
Consequently,%
\begin{equation*}
\sum\limits_{k=0}^{n}q_{n-k}\gamma _{k}=0\left( n>0\right) .
\end{equation*}%
Then
\begin{eqnarray*}
\frac{A_{n}^{\alpha _{n}}}{Q_{k}}t_{k,n}
&=&\sum\limits_{l=0}^{n-k}A_{n-k-l}^{\alpha _{n}-1}\gamma _{l} \\
&=&A_{n-k}^{\alpha _{n}-1}\gamma
_{0}+\sum\limits_{l=1}^{n-k}A_{n-k-l}^{\alpha _{n}-1}\gamma _{l} \\
&=&A_{n-k}^{\alpha _{n}-1}\left( 1+\sum\limits_{l=1}^{n-k}\frac{%
A_{n-k-l}^{\alpha _{n}-1}}{A_{n-k}^{\alpha _{n}-1}}\gamma _{l}\right)  \\
&\geq &A_{n-k}^{\alpha _{n}-1}\left( 1+\sum\limits_{l=1}^{n-k}\frac{q_{n-k-l}%
}{q_{n-k}}\gamma _{l}\right)  \\
&=&\frac{A_{n-k}^{\alpha _{n}-1}}{q_{n-k}}\left(
q_{n-k}+\sum\limits_{l=1}^{n-k}q_{n-k-l}\gamma _{l}\right)  \\
&=&\frac{A_{n-k}^{\alpha _{n}-1}}{q_{n-k}}\sum\limits_{l=0}^{n-k}q_{n-k-l}%
\gamma _{l}=0.
\end{eqnarray*}%
Hence,%
\begin{equation}
\sum\limits_{k=0}^{n}\left\vert t_{k,n}\right\vert
=\sum\limits_{k=0}^{n}t_{k,n}=1  \label{1}
\end{equation}%
and we arrive at the condition $\left( c\right) $.

Now, let us choose the sequnce $(q_n)$ as follows:
$$
q_n=\frac{1}{n+1}
$$
Then all the required conditions for $(q_n)$ are fulfilled. Moreover,
the corresponding N\"{o}rlund means coincides with the Logarithmic means. Hence, from
\begin{equation}
\sigma _{n}^{\alpha _{n}}(\mathbf{S})=\sum\limits_{k=0}^{n}t_{k,n}L_{k}%
\left( \mathbf{S}\right)   \label{2}
\end{equation}%
we immediately arrive at the required assertion.
\end{proof}

The connection between divergence of the given sequence and the divergence
of the sequence $\sigma _{n}^{\alpha _{n}}(\mathbf{S})$ was established by
Tetunashvili \cite{tetuna}. In particularly, the following theorem was proved.

\begin{Te}
Let  $\left\{ \alpha _{n}\right\} $ be a sequences such that
\begin{equation*}
0<\alpha _{n}\leq \frac{c}{\ln n},\quad \text{ where }0<c<\ln 2\text{ and }n>m.
\end{equation*}%
for some $m>0$.
Then for any sequence $\mathbf{S}:=\left\{ s_{n}:n\in \mathbb{N}\right\} $:
\begin{equation*}
\sup\limits_{n}\left\vert s_{n}\right\vert =+\infty ,
\end{equation*}%
implies
\begin{equation*}
\sup\limits_{n}\left\vert \sigma _{n}^{a_{n}}\right\vert =+\infty.
\end{equation*}%
\end{Te}

From the previous Theorem 5 and Theorem T we obtain next result if one takes $\alpha _{n}=\frac{1}{n+1}$.

\begin{theorem}
	\label{T4}Let $\mathbf{S}:=\left\{ s_{n}:n\in \mathbb{N}\right\} $ be a
	sequence of real numbers such that
	\begin{equation*}
	\sup\limits_{n}\left\vert s_{n}\right\vert =+\infty .
	\end{equation*}%
	Then
	\begin{equation*}
	\sup\limits_{n}\left\vert L_{n}\left( \mathbf{S}\right) \right\vert =+\infty
	.
	\end{equation*}
\end{theorem}

\section*{Conflicts of Interest} The authors declare that they have no conflicts of interest.

\section*{Declarations}

\subsection*{Ethical Approval}

Not applicable

\subsection*{Authorship Contribution} All authors contributed equally.


\begin{thebibliography}{10}
	
	\bibitem{AkhZviad2019}
	T.~Akhobadze and Sh. Zviadadze.
	\newblock A note on the generalized {C}es\'{a}ro means of trigonometric
	{F}ourier series.
	\newblock {\em Izv. Nats. Akad. Nauk Armenii Mat.}, 54(5):3--10, 2019.
	
	\bibitem{bochkarev1975fourier}
	S.~Bochkarev.
	\newblock A {F}ourier series in an arbitrary bounded orthonormal system that
	diverges on a set of positive measure.
	\newblock {\em Matematicheskii Sbornik}, 140(3):436--449, 1975.
	
	\bibitem{Bochk}
	S.~V. Bochkarev.
	\newblock Everywhere divergent {F}ourier series in the {W}alsh system and in
	multiplicative systems.
	\newblock {\em Uspekhi Mat. Nauk}, 59(1(355)):103--124, 2004.
	
	\bibitem{Buga}
	V.~M. Bugadze.
	\newblock Divergence of {F}ourier-{W}alsh series of bounded functions on sets
	of measure zero.
	\newblock {\em Mat. Sb.}, 185(7):119--127, 1994.
	
	\bibitem{Buzdalin}
	V.~V. Buzdalin.
	\newblock Unboundedly diverging trigonometric {F}ourier series of continuous
	functions.
	\newblock {\em Mat. Zametki}, 7:7--18, 1970.
	
	\bibitem{di2014lacunary}
	F.~Di~Plinio.
	\newblock Lacunary fourier and walsh--fourier series near l\^{} 1 l 1.
	\newblock {\em Collectanea mathematica}, 65:219--232, 2014.
	
	\bibitem{DoLacey2012}
	Y.~Q. Do and M.~Lacey.
	\newblock On the convergence of lacunacy {W}alsh-{F}ourier series.
	\newblock {\em Bull. Lond. Math. Soc.}, 44(2):241--254, 2012.
	
	\bibitem{gatJAT2010}
	G.~G\'{a}t.
	\newblock Almost everywhere convergence of {F}ej\'{e}r and logarithmic means of
	subsequences of partial sums of the {W}alsh-{F}ourier series of integrable
	functions.
	\newblock {\em J. Approx. Theory}, 162(4):687--708, 2010.
	
	\bibitem{gatConstr2019}
	G.~G\'{a}t.
	\newblock Ces\`aro means of subsequences of partial sums of trigonometric
	{F}ourier series.
	\newblock {\em Constr. Approx.}, 49(1):59--101, 2019.
	
	\bibitem{Gat2023}
	G.~G\'at.
	\newblock Almost everywhere divergence of ces\'aro means of subsequences of
	partial sums of trigonometric fourier series.
	\newblock {\em Mathematische Annalen}, 2023.
	
	\bibitem{GG2006-log}
	G.~G\'{a}t and U.~Goginava.
	\newblock Uniform and {$L$}-convergence of logarithmic means of
	{W}alsh-{F}ourier series.
	\newblock {\em Acta Math. Sin. (Engl. Ser.)}, 22(2):497--506, 2006.
	
	\bibitem{GatGogiAMSDiv2}
	G.~G\'{a}t and U.~Goginava.
	\newblock On the divergence of {N}\"{o}rlund logarithmic means of
	{W}alsh-{F}ourier series.
	\newblock {\em Acta Math. Sin. (Engl. Ser.)}, 25(6):903--916, 2009.
	
	\bibitem{gogimeasurezero}
	U.~Goginava.
	\newblock On the divergence of {W}alsh-{F}ej\'{e}r means of bounded functions
	on sets of measure zero.
	\newblock {\em Acta Math. Hungar.}, 121(4):359--369, 2008.
	
	\bibitem{goginava2021divergence}
	U.~Goginava and G.~Oniani.
	\newblock On the divergence of subsequences of partial {W}alsh-{F}ourier sums.
	\newblock {\em Journal of Mathematical Analysis and Applications},
	497(2):124900, 2021.
	
	\bibitem{gosselin1958divergence}
	R.~Gosselin.
	\newblock On the divergence of {F}ourier series.
	\newblock In {\em Proc. Amer. Math. Soc}, volume~9, pages 278--282, 1958.
	
	\bibitem{hardy2000divergent}
	G.~H. Hardy.
	\newblock {\em Divergent series}, volume 334.
	\newblock American Mathematical Soc., 2000.
	
	\bibitem{kahankatznelson1966ensembles}
	J.~Kahane and Y.~Katznelson.
	\newblock Sur les coefficients des s\'{e}ries de {F}ourier dont les sommes
	partielles sont positives sur un ensemble.
	\newblock {\em Studia Math.}, 44:555--562, 1972.
	
	\bibitem{Kap}
	I.~B. Kaplan.
	\newblock Ces\`aro means of variable order.
	\newblock {\em Izv. Vys\v{s}. U\v{c}ebn. Zaved. Matematika}, 1960(5
	(18)):62--73, 1960.
	
	\bibitem{KaragulyanMeasureZero}
	G.~A. Karagulyan.
	\newblock Divergence of general operators on sets of measure zero.
	\newblock {\em Colloq. Math.}, 121(1):113--119, 2010.
	
	\bibitem{katznelson1966ensembles}
	Y.~Katznelson.
	\newblock Sur les ensembles de divergence des s{\'e}ries trigonom{\'e}triques.
	\newblock {\em Studia Mathematica}, 26(3):305--306, 1966.
	
	\bibitem{kazaryan1991problem}
	K.~Kazaryan and P.~Ul'yanov.
	\newblock On a problem of alexits with respect to divergence of orthogonal
	fourier series.
	\newblock In {\em Dokl. Akad. Nauk SSSR}, volume 316, pages 542--546, 1991.
	
	\bibitem{kolmogoroff1923serie}
	A.~Kolmogoroff.
	\newblock Une s{\'e}rie de {F}ourier-{L}ebesgue divergente presque partout.
	\newblock {\em Fundamenta mathematicae}, 4(1):324--328, 1923.
	
	\bibitem{kolmogororf1926s}
	A.~Kolmogoroff.
	\newblock Une s{\'e}rie de {F}ourier-{L}ebesgue divergente partout.
	\newblock {\em CR Acad. ScL Paris}, 183, 1926.
	
	\bibitem{konyagin2006almost}
	S.~Konyagin.
	\newblock Almost everywhere convergence and divergence of {F}ourier series.
	\newblock In {\em International congress of mathematicians}, volume~2, pages
	1393--1403. Citeseer, 2006.
	
	\bibitem{kon}
	S.~V. Konyagin.
	\newblock On a subsequence of {F}ourier-{W}alsh partial sums.
	\newblock {\em Mat. Zametki}, 54(4):69--75, 158, 1993.
	
	\bibitem{Kon2005}
	S.~V. Konyagin.
	\newblock Divergence everywhere of subsequences of partial sums of
	trigonometric {F}ourier series.
	\newblock {\em Proc. Steklov Inst. Math.}, (Function Theory, suppl.
	2):S167--S175, 2005.
	
	\bibitem{Kon2006}
	S.~V. Konyagin.
	\newblock Almost everywhere convergence and divergence of {F}ourier series.
	\newblock In {\em International {C}ongress of {M}athematicians. {V}ol. {II}},
	pages 1393--1403. Eur. Math. Soc., Z\"{u}rich, 2006.
	
	\bibitem{Kon2011}
	S.~V. Konyagin.
	\newblock Almost everywhere divergence of lacunary subsequences of partial sums
	of {F}ourier series.
	\newblock {\em Proc. Steklov Inst. Math.}, 273(suppl. 1):S99--S106, 2011.
	
	\bibitem{lie2012pointwise}
	V.~Lie.
	\newblock On the pointwise convergence of the sequence of partial fourier sums
	along lacunary subsequences.
	\newblock {\em Journal of Functional Analysis}, 263(11):3391--3411, 2012.
	
	\bibitem{MorSiddJAT}
	F.~M\'{o}ricz and A.~H. Siddiqi.
	\newblock Approximation by {N}\"{o}rlund means of {W}alsh-{F}ourier series.
	\newblock {\em J. Approx. Theory}, 70(3):375--389, 1992.
	
	\bibitem{stechkin1951convergence}
	S.~Stechkin.
	\newblock On convergence and divergence of trigonometric series.
	\newblock {\em Uspekhi Mat. Nauk}, 6(2):148--149, 1951.
	
	\bibitem{stein1961limits}
	E.~Stein.
	\newblock On limits of sequences of operators.
	\newblock {\em Annals of Mathematics}, pages 140--170, 1961.
	
	\bibitem{tetuna}
	Sh. Tetunashvili.
	\newblock On divergence of {F}ourier series by some methods of summability.
	\newblock {\em J. Funct. Spaces Appl.}, pages Art. ID 542607, 9, 2012.
	
	\bibitem{totik1982divergence}
	V.~Totik.
	\newblock On the divergence of fourier series.
	\newblock {\em Publ. Math. Debrecen}, 29(3-4):251--264, 1982.
	
	\bibitem{zygmund2002trigonometric}
	A.~Zygmund.
	\newblock {\em Trigonometric series}, volume~2.
	\newblock Cambridge university press, 2002.
	
\end{thebibliography}
\end{document}